\documentclass{amsart}
\title{Solution to a BCC 2022 problem}
\author{Henry (Maya) Robert Thackeray}
\address{Department of Mathematics and Applied Mathematics, University of Pretoria, Pretoria, 0002 South Africa, maya.thackeray@up.ac.za, mayart314@outlook.com}
\begin{document}
\begin{abstract}
For positive integers $n$ and $k$ such that $k$ is at most $n$, we find an explicit one-to-one correspondence between the following two sets: the set of words consisting of $k$ $R$s, $k$ $U$s, and $n - k$ $D$s, where the first letter of the word is not $D$; and the set of subgraphs $H$ of a cycle of length $2n$ (where that cycle has differently labelled vertices) such that $H$ has $n$ edges and $k$ connected components. This solves a problem of Thomas Selig from the 29th British Combinatorial Conference held at Lancaster University in July 2022.
\end{abstract}
\maketitle
\section{Introduction}
This note solves a problem posed by Thomas Selig at the 29th British Combinatorial Conference (BCC 2022) held at Lancaster University in July 2022 \cite{LU22}.

Let the graph $C$ be the cycle of length $2n$ where the vertices of $C$ are the elements $0$, $1$, \ldots, $2n - 1$ of $\mathbb{Z}/(2n\mathbb{Z})$, and the edges of $C$ are the edges $\{i,i + 1\}$ with $i \in \mathbb{Z}/(2n\mathbb{Z})$. For each $i$ and $j$ in $\mathbb{Z}/(2n\mathbb{Z})$, define the path $P(i,j)$ in $C$ as follows: if $i = j$ then $P(i,j)$ is the no-edge path in which the only vertex is $i$, and if $i \neq j$ then $P(i,j)$ is the path $(i,i + 1,\ldots,j - 1,j)$ from $i$ to $j$. (For example, $P(-1,1)$ is the two-edge path $(2n - 1,0,1)$, since $-1 = 2n - 1$ in $\mathbb{Z}/(2n\mathbb{Z})$.)

A statement of the problem of Selig is as follows: For integers $n$ and $k$ with $1 \leq k \leq n$, find an explicit one-to-one corrrespondence between
\begin{itemize}
\item The set $W$ of words consisting of $k$ letters $R$, $k$ letters $U$, and $n - k$ letters $D$, in any order such that the first letter of the word is not $D$; and
\item The set $G$ of subgraphs of $C$ with $n$ edges and $k$ connected components.
\end{itemize}
(Instead of $W$, the problem originally referred to the set of all lattice paths in the $xy$ plane from $(0,0)$ to $(n,n)$ consisting of the following steps in any order, where the first step is not a $D$ step: $k$ $R$ steps, where each step moves one unit to the right; $k$ $U$ steps, where each step moves one unit up; and $n - k$ $D$ steps, where each step moves along a $45$-degree diagonal one unit to the right and one unit up. There is a clear explicit one-to-one correspondence between such paths and words in $W$: the letters of the word are the letters of the steps in the path, in the same order.)

We take a partition $W = W_{*} \cup W_{0}$ and a partition $G = G_{*} \cup G_{0}$. We explicitly describe a bijection from $W_{*}$ to $G_{*}$, and then a bijection from $W_{0}$ to $G_{0}$. (Checking that these maps are bijections is left as a straightforward exercise for the reader.) Combining these bijections gives an explicit bijection from $W$ to $G$.

\section{Partitions and sequences of numbers}

Let $W_{*}$ consist of the words $w$ in $W$ such that for some letter $D$ in $w$, there are no $R$s before that $D$ in $w$ (in such a word, the first letter is $U$); take the partition $W = W_{*} \cup W_{0}$, where $W_{0} = W - W_{*}$. Let $G_{*}$ consist of the subgraphs of $C$ that do not contain the vertex $0$; take the partition $G = G_{*} \cup G_{0}$, where $G_{0} = G - G_{*}$.

For each word $w$ in $W$, define integers $p_{0}$, \ldots, $p_{k}$ and $q_{0}$, \ldots, $q_{k}$ as follows.
\begin{itemize}
\item The number $p_{0}$ is the number of non-$U$ letters before the first $U$ in $w$.
\item The number $q_{0}$ is $1$ less than the number of $D$s before the first $R$ in $w$.
\item For positive integers $i \leq k - 1$, $p_{i}$ (respectively, $q_{i}$) is the number of non-$U$ letters (resp.\@ $D$s) between the $i$th and $(i + 1)$st $U$s (resp.\@ $R$s) in $w$.
\item The number $p_{k}$ (respectively, $q_{k}$) is the number of non-$U$ letters (resp.\@ $D$s) after the last $U$ (resp.\@ $R$) in $w$.
\end{itemize}

\section{The first bijection}

The set $W_{*}$ has $\binom{n + k - 1}{k - 1}\binom{n - 1}{k}$ words. (Proof: For each word in $W_{*}$, the first letter is $U$ and the first non-$U$ letter is $D$. To specify a word in $W_{*}$, choose $k - 1$ positions in the word for the other $U$s, then choose $k$ positions for the $R$s. We take $\binom{n - 1}{n} = 0$.) For each $w$ in $W_{*}$, we have the sequence $(p_{1},\ldots,p_{k})$ of $k$ nonnegative integers that sum to $n$, and the sequence $(q_{0},\ldots,q_{k})$ of $k + 1$ nonnegative integers that sum to $n - k - 1$, defined as above. In the element of $G_{*}$ corresponding to $w$, the connected components are the following $k$ paths, where $m \in \{1,\ldots,k\}$:
\[\begin{array}{c}
P\left(m + \sum_{i = 1}^{m - 1}p_{i} + \sum_{i = 0}^{m - 1}q_{i},m + \sum_{i = 1}^{m}p_{i} + \sum_{i = 0}^{m - 1}q_{i}\right).
\end{array}\]
For $m \in \{1,\ldots,k\}$, the $m$th path has $p_{m}$ edges; for $m \in \{1,\ldots,k - 1\}$, there are $q_{m} + 1$ edges between the $m$th and $(m + 1)$st paths. ($\sum_{i = 1}^{0}$ is the empty sum $0$.)

\section{The second bijection}

The set $W_{0}$ has $\binom{n + k}{k}\binom{n - 1}{k - 1}$ words. (Proof: The first non-$U$ letter is $R$. Choose $k$ positions for the $U$s; choose $k - 1$ positions for the other $R$s.) For each $w$ in $W_{0}$, we have the sequence $(p_{0},\ldots,p_{k})$ of $k + 1$ nonnegative integers that sum to $n$, and the sequence $(q_{1},\ldots,q_{k})$ of $k$ nonnegative integers that sum to $n - k$, defined as above. In the element of $G_{0}$ corresponding to $w$, the connected components are:
\[\begin{array}{c}
P\left(m + \sum_{i = 0}^{m - 1}p_{i} + \sum_{i = 1}^{m}q_{i},m + \sum_{i = 0}^{m}p_{i} + \sum_{i = 1}^{m}q_{i}\right)
\end{array}\]
where $m \in \{1,\ldots,k - 1\}$, and $P\left(-p_{k},p_{0}\right)$. For $m \in \{1,\ldots,k - 1\}$, the $m$th path has $p_{m}$ edges and there are $q_{m + 1} + 1$ edges between the $m$th and $(m + 1)$st paths; the path $P(-p_{k},p_{0})$ has $p_{k} + p_{0}$ edges.

\section*{Acknowledgements}

At BCC 2022, I presented postdoctoral research that I carried out at the University of Pretoria. Many thanks to everyone at the University of Pretoria, including my postdoctoral supervisor James Raftery, Roumen Anguelov, Jan Harm van der Walt, Mapundi Banda, and Anton Str\"{o}h, for their past and continuing generous support. Many thanks to Tony Nixon, Sean Prendiville, and everyone else behind the scenes of BCC 2022 for making this wonderful conference possible. Many thanks to Thomas Selig for posing this paper's problem.

\end{document}